\documentclass[12pt]{article}
\usepackage{amssymb,amsmath,amsthm,amsfonts,enumerate, bbm, mathdots}

\usepackage{latexsym}
\usepackage{amscd}
\usepackage{stmaryrd}
\usepackage{color}
\usepackage[all]{xypic}
\usepackage{epsfig}
\usepackage{graphics}
\usepackage{ifthen}
\usepackage{varioref}
\usepackage{rotating}
\usepackage{extarrows}
\usepackage{cite}
\usepackage{mathrsfs}

\numberwithin{equation}{section}

\theoremstyle{plain}
\newtheorem{thm}{Theorem}[section]
\newtheorem{cor}[thm]{Corollary}
\newtheorem{lem}[thm]{Lemma}

\newtheorem{conj}{Conjecture}[section]

\makeatletter
\newenvironment{proofof}[1]{\par
  \pushQED{\qed}%
  \normalfont \topsep6\p@\@plus6\p@\relax
  \trivlist
  \item[\hskip\labelsep
        \bfseries
    Proof of #1\@addpunct{.}]\ignorespaces
}{%
  \popQED\endtrivlist\@endpefalse
}
\makeatother

\definecolor{darkgreen}{rgb}{0.0625,0.64,0.0625}
\usepackage[
            pdfstartview=FitH,
            CJKbookmarks=true,
            bookmarksnumbered=true,
            bookmarksopen=true,
            colorlinks,
            pdfborder=001,
            linkcolor=blue,
            anchorcolor=green,
            citecolor=red
            ]{hyperref}

\newfont{\scyr}{wncyr10 scaled 550}

\def\proof{\noindent {\bf Proof.\;}}

\def\wt{\operatorname{wt}}
\def\dep{\operatorname{dep}}
\def\height{\operatorname{ht}}

\allowdisplaybreaks

\begin{document}

\title{A note on the sum of finite multiple harmonic $q$-series on $r\text{-}(r+1)$ indices}

\date{\small ~ \qquad\qquad School of Mathematical Sciences, Tongji University \newline No. 1239 Siping Road,
Shanghai 200092, China}

\author{Zhonghua Li\thanks{E-mail address: zhonghua\_li@tongji.edu.cn} ~and ~Zhenlu Wang\thanks{E-mail address: zlw@tongji.edu.cn}}

\maketitle

\begin{abstract}
We study the sum of the finite multiple harmonic $q$-series on $r\text{-}(r+1)$ indices at roots of unity with $r=1,2,3$. And we give the equivalent conditions of two conjectures regarding cyclic sums of finite multiple harmonic $q$-series on $1\text{-}2\text{-}3$ indices at roots of unity, posed recently by Kh. Pilehrood, T. Pilehrood and R. Tauraso.
\end{abstract}

{\small

}


\section{Introduction}\label{Sec:Intro}

Let $n\in \mathbb{N}$ be fixed, where $\mathbb{N}$ is the set of positive integers.  Let $q$ be a complex number satisfying $q^m\neq 1$ for any $m\in \mathbb{N}$ with $m<n$. For any multi-index $\mathbf{k}=(k_1,\ldots,k_d)\in \mathbb{N}^d$ with $d\in \mathbb{N}$,  the finite multiple harmonic $q$-series are defined by
$$z_n(\mathbf{k};q)=z_n(k_1,\ldots,k_d;q)=\sum\limits_{n>m_1>\cdots>m_d>0}\frac{q^{(k_1-1)m_1+\cdots+(k_d-1)m_d}}{[m_1]^{k_1}\cdots[m_d]^{k_d}}.$$
Here for any $m\in \mathbb{N}$, $[m]$ is the $q$-integer $[m]=\frac{1-q^m}{1-q}$. H. Bachmann, Y. Takeyama and K. Tasaka studied the special values of the finite multiple harmonic $q$-series at roots of unity in \cite{B-T-T,B-T-T-2}. Z. Li and E. Pan  introduced the interpolated finite multiple harmonic $q$-series and considered the generating function of the sums of the interpolated finite multiple harmonic $q$-series with fixed weight, depth and $i$-height in \cite{Li-Pan}. Based on the main result of \cite{Li-Pan}, we study the sum of finite multiple harmonic $q$-series on $r\text{-}(r+1)$ indices at root of unity in this paper.

Now let $\zeta_n$ be a fixed primitive $n$-th root of unity. Let $r\in \mathbb{N}$ and $l, s\in\mathbb{Z}_{\geq 0}$, where $\mathbb{Z}_{\geq 0}$ is the set of  the non-negative integers. We set
$$W_n^r(l,s)=\sum\limits_{d_0+d_1+\cdots+d_s=l\atop d_0,\ldots,d_s\geq 0}z_n\left(\{r\}^{d_0},r+1,\{r\}^{d_{1}},r+1,\ldots,r+1,\{r\}^{d_{s}};\zeta_n\right),$$
which is the sum of the finite multiple harmonic $q$-series on $r\text{-}(r+1)$ indices at roots of unity.  For example, when $r=1,2,3$, we have
\begin{align*}
W_n^1(l,s)=&\sum\limits_{d_0+d_1+\cdots+d_s=l\atop d_0,\ldots,d_s\geq 0}z_n\left(\{1\}^{d_0},2,\{1\}^{d_{1}},2,\ldots,2,\{1\}^{d_{s}};\zeta_n\right),\\
W_n^2(l,s)=&\sum\limits_{d_0+d_1+\cdots+d_s=l\atop d_0,\ldots,d_s\geq 0}z_n\left(\{2\}^{d_0},3,\{2\}^{d_{1}},3,\ldots,3,\{2\}^{d_{s}};\zeta_n\right),\\
W_n^3(l,s)=&\sum\limits_{d_0+d_1+\cdots+d_s=l\atop d_0,\ldots,d_s\geq 0}z_n\left(\{3\}^{d_0},4,\{3\}^{d_{1}},4,\ldots,4,\{3\}^{d_{s}};\zeta_n\right).
\end{align*}

We also give the following notation. Let
\begin{align*}
A(n,l,s)&=\sum\limits_{i=1}^{s+l+1}\frac{-(s+1)}{i\binom{s+l}{l}}\binom{i}{s+l+1-i}\sum\limits_{j=0}^{2i-s-l-1}(-4)^{j}\binom{2i-s-l-1}{j}\\
&\qquad\qquad\qquad \times\binom{n+2s+2l+1-i+j}{4s+4l+3-2i+2j}\binom{2s+2l+2-2i}{s+1-j}.
\end{align*}
One can see that the right-hand side of the above formula is complicated. It is better to give a concise one, and we are supposed to find a suitable solution.

The main results of this paper are the following theorems.

\begin{thm}\label{Theorem:W_n^1 sum formula}
For $n\in \mathbb{N}$ and $l,s\in \mathbb{Z}_{\geq 0}$, we have
  \begin{align*}
  W_n^1(l,s)=\frac{(-1)^s}{n(s+1)}\binom{s+l}{l}\binom{n+s}{2s+l+1}(1-\zeta_n)^{2s+l}.
  \end{align*}
\end{thm}

\begin{thm}\label{Theorem:W_n^2 sum formula}
For $n\in \mathbb{N}$ and $l,s\in \mathbb{Z}_{\geq 0}$, we have
  \begin{align*}
  W_n^2(l,s)=\frac{(-1)^l}{n^2(s+1)}\binom{s+l}{l}\left[\binom{n+2s+l+1}{3s+2l+2}+(-1)^s\binom{n+s+l}{3s+2l+2}\right](1-\zeta_n)^{3s+2l}.
  \end{align*}
\end{thm}

\begin{thm}\label{Theorem:W_n^3 sum formula}
For $n\in \mathbb{N}$ and $l,s\in \mathbb{Z}_{\geq 0}$, we have
  \begin{align*}
  W_n^3(l,s)=&\frac{(-1)^s}{n^3(s+1)}\binom{s+l}{l}\left[\binom{n+3s+2l+2}{4s+3l+3}\right.\\
  &\qquad\qquad\left.+(-1)^l\binom{n+s+l}{4s+3l+3}+A(n,l,s)\right](1-\zeta_n)^{4s+3l}.
  \end{align*}
\end{thm}

The above theorems are related to a recent work \cite{P-P-Tauraso} of Kh. Pilehrood, T. Pilehrood and R. Tauraso. In \cite{P-P-Tauraso}, the cyclic sums of finite multiple harmonic $q$-series on $1\text{-}2\text{-}3$ indices at roots of unity are studied and two conjectures are proposed. For $r\in\mathbb{N}$ and $s, d_0, d_1,\ldots,d_s\in \mathbb{Z}_{\geq 0}$, we define the following cyclic sum  $$C_n^r(d_0,d_1,\ldots,d_s)=\sum\limits_{j=0}^{s}z_n\left(\{r\}^{d_j},r+1,\{r\}^{d_{j+1}},r+1,\ldots,r+1,\{r\}^{d_{j+s}};\zeta_n\right),$$
where $d_j=d_i$ if $j\equiv i$ modulo $s+1$. The two conjectures given in \cite{P-P-Tauraso} are

\begin{conj}\label{Conjecture:r=1}
For $s, d_0, d_1,\ldots,d_s\in \mathbb{Z}_{\geq 0}$, set $k=\sum\limits_{j=0}^{s}d_j+2s$. Then for any positive integer $n$ with $n>k$, we have
\begin{align*}
C_n^1(d_0,d_1,\ldots,d_s)=\frac{(-1)^s}{n}\binom{n+s}{k+1}(1-\zeta_n)^k.
\end{align*}
\end{conj}

\begin{conj}\label{Conjecture:r=2}
For $s, d_0, d_1,\ldots,d_s\in \mathbb{Z}_{\geq 0}$, set $k=2\sum\limits_{j=0}^{s}d_j+3s$. Then for any positive integer $n$ with $n>k$, we have
\begin{align*}
C_n^2(d_0,d_1,\ldots,d_s)\in (1-\zeta_n)^k\mathbb{Q}.
\end{align*}
\end{conj}

It is obvious that Conjecture \ref{Conjecture:r=1} implies that the cyclic sum $C_n^1(d_0,d_1,\ldots,d_s)$ depends only on the number and the sum of its
arguments. Using Theorem \ref{Theorem:W_n^1 sum formula}, we prove that the converse of this fact is also true.

\begin{cor}\label{Corollary:equivalent conj r=1}
For $n\in\mathbb{N}$ and $k,s\in\mathbb{Z}_{\geq 0}$ with $n>k$,  the following two statements are equivalent:
\begin{description}
  \item[(i)] for any $d_0, d_1,\ldots,d_s\in \mathbb{Z}_{\geq 0}$ with $\sum\limits_{j=0}^{s}d_j=k-2s$, it holds
  \begin{align}\label{Eq:conj C^n_1}
C_n^1(d_0,d_1,\ldots,d_s)=\frac{(-1)^s}{n}\binom{n+s}{k+1}(1-\zeta_n)^k;
\end{align}
  \item[(ii)] for any $d_0, d_1,\ldots,d_s,d_0',d_1',\ldots,d_s'\in \mathbb{Z}_{\geq 0}$ with $\sum\limits_{j=0}^{s}d_j=\sum\limits_{j=0}^{s}d_j'=k-2s$, it holds
  $$C_n^1(d_0,d_1,\ldots,d_s)=C_n^1(d_0',d_1',\ldots,d_s').$$
\end{description}
\end{cor}

Similarly, from Theorem \ref{Theorem:W_n^2 sum formula} and Theorem \ref{Theorem:W_n^3 sum formula}, we get

\begin{cor}\label{Corollary:equivalent conj r=2}
For $n\in\mathbb{N}$ and $k,s\in\mathbb{Z}_{\geq 0}$ with $n>k$,  the following two statements are equivalent:
\begin{description}
  \item[(i)] for any $d_0, d_1,\ldots,d_s\in \mathbb{Z}_{\geq 0}$ with $2\sum\limits_{j=0}^{s}d_j=k-3s$, it holds
  \begin{align*}
C_n^2(d_0,d_1,\ldots,d_s)=\frac{(-1)^{d-s}}{n^2}\left[\binom{n+k-d+1}{k+2}+(-1)^s\binom{n+d}{k+2}\right](1-\zeta_n)^k,
\end{align*}
where $d=\sum\limits_{j=0}^sd_j+s$;
  \item[(ii)] for any $d_0, d_1,\ldots,d_s,d_0',d_1',\ldots,d_s'\in \mathbb{Z}_{\geq 0}$ with $2\sum\limits_{j=0}^{s}d_j=2\sum\limits_{j=0}^{s}d_j'=k-3s$, it holds
  $$C_n^2(d_0,d_1,\ldots,d_s)=C_n^2(d_0',d_1',\ldots,d_s').$$
\end{description}
\end{cor}

\begin{cor}\label{Corollary:equivalent conj r=3}
For $n\in\mathbb{N}$ and $k,s\in\mathbb{Z}_{\geq 0}$ with $n>k$,  the following two statements are equivalent:
\begin{description}
  \item[(i)] for any $d_0, d_1,\ldots,d_s\in \mathbb{Z}_{\geq 0}$ with $3\sum\limits_{j=0}^{s}d_j=k-4s$, it holds
  \begin{align*}
C_n^3(d_0,d_1,\ldots,d_s)=&\frac{(-1)^s}{n^3}\left[\binom{n+k-d+2}{k+3}\right.\\
&\quad\left.+(-1)^{d-s}\binom{n+d}{k+3}+A(n,d-s,s)\right](1-\zeta_n)^k,
\end{align*}
where $d=\sum\limits_{j=0}^sd_j+s$;
  \item[(ii)] for any $d_0, d_1,\ldots,d_s,d_0',d_1',\ldots,d_s'\in \mathbb{Z}_{\geq 0}$ with $3\sum\limits_{j=0}^{s}d_j=3\sum\limits_{j=0}^{s}d_j'=k-4s$, it holds
  $$C_n^3(d_0,d_1,\ldots,d_s)=C_n^3(d_0',d_1',\ldots,d_s').$$
\end{description}
\end{cor}

Based on Corollary \ref{Corollary:equivalent conj r=2}, it is natural to guess that the explicit expression of $C_n^2$ in Conjecture \ref{Conjecture:r=2} should be
\begin{align*}
C_n^2(d_0,d_1,\ldots,d_s)\overset{?}{=}\frac{(-1)^{d-s}}{n^2}\left[\binom{n+k-d+1}{k+2}+(-1)^s\binom{n+d}{k+2}\right](1-\zeta_n)^k,\\
\end{align*}
where $k=2\sum\limits_{j=0}^sd_j+3s$ and $d=\sum\limits_{j=0}^sd_j+s$. And from Corollary \ref{Corollary:equivalent conj r=3}, we may guess that
\begin{align*}
C_n^3(d_0,d_1,\ldots,d_s)\overset{?}{=}&\frac{(-1)^s}{n^3}\left[\binom{n+k-d+2}{k+3}\right.\\
&\qquad\left.+(-1)^{d-s}\binom{n+d}{k+3}+A(n,d-s,s)\right](1-\zeta_n)^k,
\end{align*}
where $k=3\sum\limits_{j=0}^sd_j+4s$ and $d=\sum\limits_{j=0}^sd_j+s$.

In Section \ref{Sec:generating function}, we study the generating function of finite multiple harmonic $q$-series on $r\text{-}(r+1)$ indices at roots of unity with fixed weight, depth and $i$-height. Several special cases are discussed in Section \ref{Sec:special cases}. In Subsection \ref{Subsec:r=1}, we prove Theorem \ref{Theorem:W_n^1 sum formula} and Corollary \ref{Corollary:equivalent conj r=1} by considering the case of $r=1$. In Subsection \ref{Subsec:r=2}, we prove Theorem \ref{Theorem:W_n^2 sum formula} and Corollary \ref{Corollary:equivalent conj r=2} by considering the case of $r=2$. In Subsection \ref{Subsec:r=3}, we prove Theorem \ref{Theorem:W_n^3 sum formula} and Corollary \ref{Corollary:equivalent conj r=3} by considering the case of $r=3$.


\section{Generating function}\label{Sec:generating function}

In this section, we first recall \cite[Theorem 1.1]{Li-Pan}. For a multi-index $\mathbf{k}=(k_1,\ldots,k_d)\in \mathbb{N}^d$,  we define its weight and depth respectively by $$\wt(\mathbf{k})=k_1+\cdots+k_d,\quad \dep(\mathbf{k})=d.$$
And for an integer $i\in \mathbb{N}$, the $i$-height is defined by $$i\text{-}\height(\mathbf{k})=\#\{j\mid k_j\geq i+1\}.$$
Let $t$ be a formal parameter. Then the interpolated finite multiple harmonic $q$-series $z_n^{t}(\mathbf{k};q)$ is defined as
$$z_n^{t}(\mathbf{k};q)=z_n^{t}(k_1,\ldots,k_d;q)=\sum\limits_{\mathbf{p}}(1-q)^{k-\wt(\mathbf{p})}z_n(\mathbf{p};q)t^{d-\dep(\mathbf{p})},$$
where $k=\wt(\mathbf{k})$ and $\sum\limits_{\mathbf{p}}$ is the sum where $\mathbf{p}$ runs over all multi-indices of the form $\mathbf{p}=(k_1\square\cdots\square k_d)$, in which each $\square$ is filled by ``,",``+" or ``$-1+$". It is easy to see that
$$z_n^{0}(\mathbf{k};q)=z_n(\mathbf{k};q),\quad z_n^{1}(\mathbf{k};q)=z_n^{\star}(\mathbf{k};q).$$
Here $z_n^{\star}(\mathbf{k};q)$ is the finite multiple harmonic $q$-series of star-version defined by
$$z_n^{\star}(\mathbf{k};q)=z_n^{\star}(k_1,\ldots,k_d;q)=\sum\limits_{n>m_1\geq\cdots\geq m_d>0}\frac{q^{(k_1-1)m_1+\cdots+(k_d-1)m_d}}{[m_1]^{k_1}\cdots[m_d]^{k_d}}.$$
In \cite{Li-Pan}, instead of $z_n^{t}(\mathbf{k};q)$ a modified version was used, which is defined as
$$\bar{z}_n^t(\mathbf{k};q)=(1-q)^{-\wt(\mathbf{k})}z_n^{t}(\mathbf{k};q).$$

Let $r\in\mathbb{N}$ be fixed. In \cite{Li-Pan}, Z. Li and E. Pan considered the generating function of the sums of the interpolated finite multiple harmonic $q$-series with fixed weight, depth and $1$-height, $2$-height, $\ldots$, $r$-height. Precisely, for any non-negative integers $k,d,h_1,\ldots,h_r$, set
$$G_n^t(k,d,h_1,\ldots,h_r;q)=\sum\limits_{\mathbf{k}\in I(k,d,h_1,\ldots,h_r)}\bar{z}_n^{t}(\mathbf{k};q),$$
where $I(k,d,h_1,\ldots,h_r)$ is the set of multi-indices of weight $k$, depth $d$, $1$-height $h_1$, $\ldots$, $r$-height $h_r$. The sum is treated as $0$ whenever the index set is empty except for $G_n^t(0,0,\ldots,0;q)=1$. For formal variables $u_1,\ldots,u_{r+2}$, we define the generating function
\begin{align*}
\Psi_n^t(q)&=\Psi_n^t(u_1,\ldots,u_{r+2};q)\\
&=\sum\limits_{k,d,h_1,\ldots,h_r\geq 0}G_n^t(k,d,h_1,\ldots,h_r;q)u_1^{k-d-\sum\limits_{j=1}^{r}h_j}u_2^{d-h_1}u_3^{h_1-h_2}\cdots u_{r+1}^{h_{r-1}-h_r}u_{r+2}^{h_r}.
\end{align*}
Then the generating function was computed in \cite[Theorem 1.1]{Li-Pan}.

\begin{thm}[{\cite[Theorem 1.1]{Li-Pan}}]\label{Thm:Li-Pan}
Let $r\in \mathbb{N}$ and $u_1,\ldots,u_{r+2}$ be variables. Set
\begin{align*}
x_1=\frac{u_1}{1+u_1}
\end{align*} and
\begin{align*}
x_i=\sum\limits_{j=i}^{r+1}(-1)^{j-i}\binom{j-2}{i-2}\left(u_j-\frac{u_{r+2}}{u_1^{r+2-j}}\right)+\frac{u_{r+2}}{u_1^{r+2-i}(1+u_1)^{i-1}}
\end{align*}
for $i=2,\cdots,r+2$. Then we have
$$\Psi_n^t(q)=\frac{\prod\limits_{j=1}^{n-1}P^{t-1}(1-q^j)}{\prod\limits_{j=1}^{n-1}P^t(1-q^j)},$$ where
\begin{align}
P^t(T)=T^{r+1}-(x_1+tx_2)T^r-t\sum\limits_{i=0}^{r-1}(x_{r+2-i}-x_1x_{r+1-i})T^i.
\end{align}\label{Eq:P^t(T)}
\end{thm}

Now we consider the sums of the interpolated finite multiple harmonic $q$-series on $r\text{-}(r+1)$ indices. It is enough to set $u_1=\cdots=u_r=0$. And in this case, we have the following lemma.

\begin{lem}\label{Lemma:x_i}
If $u_1=\cdots=u_r=0$, then we have $x_1=0$ and
\begin{align*}
x_i=(-1)^{r+1-i}\binom{r-1}{i-2}u_{r+1}+(-1)^{r-i}\binom{r}{i-2}u_{r+2}
\end{align*}
for $i=2,\ldots,r+2$.
\end{lem}

\proof
We get the result by the same way as in \cite[Lemma 4.1]{Li-Pan}.
\qed

Using Lemma \ref{Lemma:x_i} and \eqref{Eq:P^t(T)}, if $u_1=\cdots=u_r=0$, we have
\begin{align*}
P^t(T)&=T^{r+1}-tx_2T^r-t\sum\limits_{i=0}^{r-1}x_{r+2-i}T^i\\
&=T^{r+1}-tT(1-T)^{r-1}u_{r+1}-t(1-T)^ru_{r+2}.
\end{align*}
Hence we get the following corollary, which is about the generating function of the sums of the interpolated finite multiple harmonic $q$-series on $r\text{-}(r+1)$ indices at roots of unity.

\begin{cor}
For $r\in \mathbb{N}$ and formal variables $u,v$, we have
$$\sum\limits_{l,s\geq 0}\sum\limits_{a_0+a_1+\cdots+a_s=l\atop a_0,\ldots,a_s\geq 0}\bar{z}_n^t\left(\{r\}^{a_0},r+1,\{r\}^{a_{1}},r+1,\ldots,r+1,\{r\}^{a_{s}};\zeta_n\right)u^lv^s=\frac{\prod\limits_{j=1}^{n-1}\widetilde{P}^{t-1}(\zeta_n^j)}{\prod\limits_{j=1}^{n-1}\widetilde{P}^t(\zeta_n^j)},$$
where $$\widetilde{P}^{t}(T)=(1-T)^{r+1}-t(1-T)T^{r-1}u-tT^rv.$$
\end{cor}

Instead of $\prod\limits_{j=1}^{n-1}\widetilde{P}^t(\zeta_n^j)$, it is convenient to use
$$H(n,r,t)=-\frac{v}{n}\prod\limits_{j=1}^{n-1}\widetilde{P}^t(\zeta_n^j).$$
Since we want to study the sum of finite multiple harmonic $q$-series at roots of unity, the following lemma is necessary.

\begin{lem}\label{Lem: H_n^t t=0}
$H(n,r,0)=-n^rv$.
\end{lem}

\proof
By the definition of $H(n,r,t)$, we find
\begin{align*}
H(n,r,0)=-\frac{v}{n}\prod_{j=1}^{n-1}(1-\zeta_n^j)^{r+1}=-\frac{v}{n}\cdot n^{r+1}=-n^rv,
\end{align*}
as desired.
\qed

Now we set $$\widetilde{P}^t(T)=(\beta_1^t-T)(\beta_2^t-T)\cdots(\beta_{r+1}^t-T).$$ Then the elementary symmetric polynomials of $\beta_1^t,\ldots,\beta_{r+1}^t$ are given by
\begin{align}\label{Eq:beta}
e_{r,j}=\sum\limits_{1\leq i_1<\cdots<i_j\leq r+1}\beta_{i_1}^t\cdots\beta_{i_j}^t=
\begin{cases}
r+1+(-1)^rt(u-v) & \text{if\;} j=1, \\
\binom{r+1}{2}+(-1)^rtu & \text{if\;} j=2, \\
\binom{r+1}{j} & \text{if\;} j=3,\cdots,r+1.
\end{cases}
\end{align}
Using the formula $\prod_{j=1}^{n-1}(T-\zeta_n^j)=\frac{T^n-1}{T-1}$, we have
\begin{align}\label{Eq:P(T)}
H(n,r,t)=&-\frac{v}{n}\prod\limits_{j=1}^{n-1}\prod\limits_{i=1}^{r+1}(\beta_i^t-\zeta_n^j)=-\frac{v}{n}\prod\limits_{i=1}^{r+1}
\frac{(\beta_i^t)^n-1}{\beta_i^t-1}\notag \\=&\frac{\prod\limits_{i=1}^{r+1}[(\beta_i^t)^n-1]}{tn}.
\end{align}
Let $w$ be a variable. We consider the generating function
\begin{align*}
&\sum\limits_{n=1}^{\infty}H(n,r,t)w^n=\sum\limits_{n=1}^{\infty}\frac{\prod\limits_{i=1}^{r+1}[(\beta_i^{t})^n-1]}{tn}w^n\\
=&\frac{1}{t}\sum\limits_{n=1}^{\infty}\frac{w^n}{n}
\left[(-1)^{r+1}+\sum\limits_{j=1}^{r+1}(-1)^{r+1-j}\sum\limits_{1\leq i_1<\cdots<i_j\leq r+1}(\beta_{i_1}^t\cdots\beta_{i_j}^t)^n\right]\\
=&\frac{(-1)^r}{t}\log(1-w)+\sum\limits_{j=1}^{r+1}\frac{(-1)^{r-j}}{t}\sum\limits_{1\leq i_1<\cdots<i_j\leq r+1}\log(1-\beta_{i_1}^t\cdots\beta_{i_j}^tw).
\end{align*}
Hence we get following result.

\begin{thm}\label{Thm:application}
For $r\in \mathbb{N}$ and formal variables $u,v$, let $\beta_1^t,\ldots,\beta_{r+1}^t$ be determined by \eqref{Eq:beta}. Let $\widetilde{F}(t,r,0)=1-w$ and
\begin{align}\label{Eq:F^t}
\widetilde{F}(t,r,j)=\prod\limits_{1\leq i_1<\cdots<i_j\leq r+1}(1-\beta_{i_1}^t\cdots\beta_{i_j}^tw),\quad (1\leq j\leq r+1).
\end{align}
Then we have
\begin{align*}
\sum\limits_{l,s\geq 0}\sum\limits_{a_0+a_1+\cdots+a_s=l\atop a_0,\ldots,a_s\geq 0}\bar{z}_n^t\left(\{r\}^{a_0},r+1,\{r\}^{a_{1}},r+1,\ldots,r+1,\{r\}^{a_{s}};\zeta_n\right)u^lv^s=\frac{H(n,r,t-1)}{H(n,r,t)},
\end{align*}
where $H(n,r,t)$ is determined by $$\sum\limits_{n=1}^{\infty}H(n,r,t)w^n=\frac{(-1)^r}{t}\log\left(\prod\limits_{j=0}^{r+1}\widetilde{F}(t,r,j)^{(-1)^j}\right).$$
And in particular, we have
$$\sum\limits_{l,s\geq 0}W_n^r(l,s)(1-\zeta_n)^{-(r+1)s-rl}u^lv^s=-\frac{H(n,r,-1)}{n^rv}.$$
\end{thm}

Therefore the rest of the task is to compute $H(n,r,-1)$. We do in the special cases $r=1,2,3$ in the next section.


\section{Special cases}\label{Sec:special cases}

\subsection{The case of $r=1$}\label{Subsec:r=1}

Set $r=1$ in Theorem \ref{Thm:application}, we can prove Theorem \ref{Theorem:W_n^1 sum formula}.

\begin{proofof}{Theorem \ref{Theorem:W_n^1 sum formula}}
Set $r=1$. By \eqref{Eq:beta}, we have $$e_{1,1}=2-t(u-v),\quad e_{1,2}=1-tu.$$
Then by \eqref{Eq:F^t}, we get
\begin{align*}
\widetilde{F}(t,1,1)=(1-w)^2+t(u-v)w-tuw^2,\quad \widetilde{F}(t,1,2)=1-w+tuw,
\end{align*}
which implies
\begin{align*}
&\sum\limits_{n=1}^{\infty}H(n,1,t)w^n=-\frac{1}{t}\log\frac{(1-w)(1-w+tuw)}{(1-w)^2+t(u-v)w-tuw^2}\\
=&-\frac{1}{t}\log\left(1-\frac{tuw^2-tuw}{(1-w)^2}\right)+\frac{1}{t}\log\left(1-\frac{tuw^2-t(u-v)w}{(1-w)^2}\right)\\
=&\sum\limits_{i=1}^{\infty}\sum\limits_{m=0}^i\sum\limits_{a=0}^{\infty}\frac{1}{i}(-1)^{i-m+1}\binom{i}{m}\binom{2i+a-1}{2i-1}t^{i-1}\left[u^m(u-v)^{i-m}-u^i\right]w^{i+m+a}.
\end{align*}
Here we have used the formula
$$(1-w)^{-2i}=\sum\limits_{a=0}^\infty\binom{2i+a-1}{2i-1}w^a.$$
Hence we have
\begin{align*}
H(n,1,t)=\sum\limits_{i=1}^n\sum\limits_{m=0}^i\frac{1}{i}(-1)^{i-m+1}\binom{i}{m}\binom{n+i-1-m}{2i-1}t^{i-1}\left[u^m(u-v)^{i-m}-u^i\right].
\end{align*}
Using Theorem \ref{Thm:application}, we get
\begin{align*}
&\sum\limits_{l,s\geq 0}W_n^1(l,s)(1-\zeta_n)^{-2s-l}u^lv^s\\
=&\sum\limits_{i=1}^n\sum\limits_{m=0}^i\frac{(-1)^{m+1}}{ni}\binom{i}{m}\binom{n+i-1-m}{2i-1}\left[u^m(u-v)^{i-m}-u^i\right]v^{-1},
\end{align*}
which implies
\begin{align}
&\sum\limits_{l,s\geq 0}W_n^1(l,s)(1-\zeta_n)^{-2s-l}u^lv^s\nonumber\\
=&\sum\limits_{i=1}^n\sum\limits_{m=0}^i\sum\limits_{a=1}^{i-m}\frac{(-1)^{a+m+1}}{ni}\binom{i}{m}\binom{n+i-1-m}{2i-1}\binom{i-m}{a}u^{i-a}v^{a-1}.\label{Eq:W_n1}
\end{align}
The coefficient of $u^lv^s$ in the right-hand side of \eqref{Eq:W_n1} is
\begin{align*}
&\sum\limits_{m=0}^l\frac{(-1)^{s+m}}{n(s+l+1)}\binom{s+l+1}{m}\binom{n+s+l-m}{2s+2l+1}\binom{s+l+1-m}{s+1}\\
=&\sum\limits_{m=0}^l\frac{(-1)^{s+m}}{n(s+l+1)}\binom{s+l+1}{l}\binom{l}{m}\binom{n+s+l-m}{2s+2l+1}\\
=&\frac{(-1)^s}{n(s+l+1)}\binom{s+l+1}{l}\sum\limits_{m=0}^l(-1)^m\binom{l}{m}\binom{n+s+l-m}{2s+2l+1}.
\end{align*}
Using the combinatorial identity \cite[(3.49)]{Gould}
\begin{align}\label{Eq:zuheshu identity}
\sum\limits_{k=0}^n(-1)^k\binom{n}{k}\binom{p-k}{q}=\binom{p-n}{q-n},
\end{align}
we find the coefficient of $u^lv^s$ in the right-hand side \eqref{Eq:W_n1} is
\begin{align*}
\frac{(-1)^s}{n(s+1)}\binom{s+l}{l}\binom{n+s}{2s+l+1}.
\end{align*}
Then the theorem is proved.
\end{proofof}

We use Theorem \ref{Theorem:W_n^1 sum formula} to prove Corollary \ref{Corollary:equivalent conj r=1}.

\begin{proofof}{Corollary \ref{Corollary:equivalent conj r=1}}
By Theorem \ref{Theorem:W_n^1 sum formula}, we have
\begin{align*}
\sum\limits_{d_0+\cdots+d_s=l\atop d_0,\ldots,d_s\geq 0}C_n^1(d_0,d_1,\ldots,d_s)&=(s+1)W_n^1(l,s)\\
&=\frac{(-1)^s}{n}\binom{s+l}{l}\binom{n+s}{2s+l+1}(1-\zeta_n)^{2s+l}.
\end{align*}
Now if (ii) is valid, then since
$$\sum\limits_{d_0+\cdots+d_s=l\atop d_0,\ldots,d_s\geq 0}1=\binom{s+l}{l},$$
we find (i) is also true.
\end{proofof}

\subsection{The case of $r=2$}\label{Subsec:r=2}

In this subsection, we consider the case of $r=2$ and prove Theorem \ref{Theorem:W_n^2 sum formula} and Corollary \ref{Corollary:equivalent conj r=2}.

\begin{proofof}{Theorem \ref{Theorem:W_n^2 sum formula}}
We have
\begin{align*}
e_{2,1}=3+t(u-v),\quad e_{2,2}=3+tu,\quad e_{2,3}=1.
\end{align*}
Then by \eqref{Eq:F^t}, we get
\begin{align*}
&\widetilde{F}(t,2,1)=(1-w)^3-t(u-v)w+tuw^2,\\
&\widetilde{F}(t,2,2)=(1-w)^3-tuw+t(u-v)w^2,\\
&\widetilde{F}(t,2,3)=1-w,
\end{align*}
which implies
\begin{align*}
&\sum\limits_{n=1}^{\infty}H(n,2,t)w^n=\frac{1}{t}\log\frac{(1-w)^3-tuw+t(u-v)w^2}{(1-w)^3-t(u-v)w+tuw^2}\\
=&\frac{1}{t}\log\left(1-\frac{tuw-t(u-v)w^2}{(1-w)^3}\right)-\frac{1}{t}\log\left(1-\frac{t(u-v)w-tuw^2}{(1-w)^3}\right)\\
=&\sum\limits_{i=1}^{\infty}\sum\limits_{m=0}^i\sum\limits_{a=0}^{\infty}
\frac{(-1)^m}{i}\binom{i}{m}\binom{3i+a-1}{3i-1}t^{i-1}\left[u^m(u-v)^{i-m}-u^{i-m}(u-v)^m\right]w^{i+m+a}.
\end{align*}
Considering the coefficient of $w^n$, we find
\begin{align*}
H(n,2,t)=\sum\limits_{i=1}^n\sum\limits_{m=0}^i\frac{(-1)^m}{i}\binom{i}{m}\binom{n+2i-1-m}{3i-1}t^{i-1}\left[u^m(u-v)^{i-m}-u^{i-m}(u-v)^m\right].
\end{align*}
Hence from Theorem \ref{Thm:application}, we get
\begin{align*}
&\sum\limits_{l,s\geq 0}W_n^2(l,s)(1-\zeta_n)^{-3s-2l}u^lv^s\\
=&\sum\limits_{i=1}^n\sum\limits_{m=0}^i\frac{(-1)^{m+i}}{n^2i}\binom{i}{m}\binom{n+2i-1-m}{3i-1}\left[u^m(u-v)^{i-m}-u^{i-m}(u-v)^m\right]v^{-1}\\
=&M_1+M_2,
\end{align*}
where
\begin{align*}
M_1=&\sum\limits_{i=1}^n\sum\limits_{m=0}^i\sum\limits_{a=1}^{i-m}
\frac{(-1)^{a+m+i}}{n^2i}\binom{i}{m}\binom{n+2i-1-m}{3i-1}\binom{i-m}{a}u^{i-a}v^{a-1},\\
M_2=&-\sum\limits_{i=1}^n\sum\limits_{m=0}^i\sum\limits_{b=1}^m
\frac{(-1)^{b+m+i}}{n^2i}\binom{i}{m}\binom{n+2i-1-m}{3i-1}\binom{m}{b}u^{i-b}v^{b-1}.
\end{align*}
The coefficient of $u^lv^s$ in $M_1$ is
\begin{align*}
\sum\limits_{m=0}^l\frac{(-1)^{m+l}}{n^2(s+l+1)}\binom{s+l+1}{l}\binom{n+2s+2l+1-m}{3s+3l+2}\binom{l}{m},
\end{align*}
and in $M_2$ is
\begin{align*}
\sum\limits_{m=0}^l\frac{(-1)^{m+s+l}}{n^2(s+l+1)}\binom{s+l+1}{l}\binom{n+s+2l-m}{3s+3l+2}\binom{l}{m}.
\end{align*}
Then using \eqref{Eq:zuheshu identity}, we finally get the result.
\end{proofof}

\begin{proofof}{Corollary \ref{Corollary:equivalent conj r=2}}
The proof is similar as that of Corollary \ref{Corollary:equivalent conj r=1}.
\end{proofof}

\subsection{The case of $r=3$}\label{Subsec:r=3}
The following Lemma is needed.
\begin{lem}\label{Lemma:combinatorial identity}
For integers $n,p,q$, we have the following combinatorial identity
\begin{align*}
\sum\limits_{k=1}^n\frac{1}{k}\binom{k}{n-k}\binom{p-k}{q-2k}=\frac{1}{n}\binom{p}{q-n}+\frac{(-1)^n}{n}\binom{p-n}{q-n}.
\end{align*}
\end{lem}
\proof
We consider the following generating function
\begin{align*}
&\sum\limits_{n=1}^{\infty}\sum\limits_{q=1}^{\infty}\sum\limits_{k=1}^n\frac{1}{k}\binom{k}{n-k}\binom{p-k}{q-2k}x^ny^q\\
=&\sum\limits_{k=1}^{\infty}\frac{1}{k}x^k(1+x)^ky^{2k}(1+y)^{p-k}\\
=&(1+y)^p\sum\limits_{k=1}^{\infty}\frac{1}{k}\left[xy^2(1+x)(1+y)^{-1}\right]^k\\
=&(1+y)^p\left[-\log(1-xy^2(1+x)(1+y)^{-1})\right].
\end{align*}
Since $xy^2(1+x)(1+y)^{-1}=(1-xy)(1+\frac{xy}{1+y})$, we have
\begin{align*}
&\sum\limits_{n=1}^{\infty}\sum\limits_{q=1}^{\infty}\sum\limits_{k=1}^n\frac{1}{k}\binom{k}{n-k}\binom{p-k}{q-2k}x^ny^q\\
=&(1+y)^p\left[-\log(1-xy)-\log(1+\frac{xy}{1+y})\right]\\
=&\sum\limits_{i=1}^{\infty}\sum\limits_{j=0}^p\frac{1}{i}\binom{p}{j}x^iy^{i+j}+\sum\limits_{i=1}^{\infty}\sum\limits_{j=0}^{p-i}\frac{(-1)^i}{i}\binom{p-i}{j}x^iy^{i+j}.
\end{align*}
Let $i=n$ and $j=q-n$, one can get the result.
\qed

Now we prove Theorem \ref{Theorem:W_n^3 sum formula} and Corollary \ref{Corollary:equivalent conj r=3} by setting $r=3$ in Theorem \ref{Thm:application}.
\begin{proofof}{Theorem \ref{Theorem:W_n^3 sum formula}}
If $r=3$, we have
\begin{align*}
e_{3,1}=4-t(u-v),\quad e_{3,2}=6-tu,\quad e_{3,3}=4,\quad e_{3,4}=1.
\end{align*}
Then by \eqref{Eq:F^t}, we get
\begin{align*}
&\widetilde{F}(t,3,1)=(1-w)^4+t(u-v)w-tuw^2,\\
&\widetilde{F}(t,3,2)=(1-w)^6+t[uw(1-w)^4+4vw^2(1-w)^2-tw^3(u-v)^2],\\
&\widetilde{F}(t,3,3)=(1-w)^4+t(u-v)w^3-tuw^2,\\
&\widetilde{F}(t,3,4)=1-w.
\end{align*}
which implies
\begin{align*}
&\sum\limits_{n=1}^{\infty}H(n,3,t)w^n\\
=&-\frac{1}{t}\log\frac{(1-w)^2\left\{(1-w)^6+t[uw(1-w)^4+4vw^2(1-w)^2-tw^3(u-v)^2]\right\}}
{[(1-w)^4+t(u-v)w-tuw^2][(1-w)^4+t(u-v)w^3-tuw^2]}\\
=&N_1^t+N_2^t+N_3^t,
\end{align*}
where
\begin{align*}
&N_1^t=-\frac{1}{t}\log\left(1+\frac{t\left[uw(1-w)^4+4vw^2(1-w)^2-tw^3(u-v)^2\right]}{(1-w)^6}\right),\\
&N_2^t=\frac{1}{t}\log\left(1-\frac{tuw^2-t(u-v)w}{(1-w)^4}\right),\\
&N_3^t=\frac{1}{t}\log\left(1-\frac{tuw^2-t(u-v)w^3}{(1-w)^4}\right).
\end{align*}
We calculate $N_1^t,N_2^t,N_3^t$ respectively. For $N_1^t$, we have
\begin{align*}
N_1^t&=\sum\limits_{i=1}^{\infty}\frac{(-1)^{i}t^{i-1}}{i(1-w)^{6i}}[uw(w-1)^4+4vw^2(w-1)^2-tw^3(u-v)^2]^i\\
&=\sum\limits_{i=1}^{\infty}\sum\limits_{a_1+a_2+a_3=i \atop a_1,a_2,a_3\geq 0}\frac{(-1)^{i+a_3}t^{i+a_3-1}4^{a_2}}{i(1-w)^{6i}}
\frac{i!}{a_1!a_2!a_3!}u^{a_1}v^{a_2}(u-v)^{2a_3}w^{a_1+2a_2+3a_3}(w-1)^{4a_1+2a_2}\\
&=\sum\limits_{i=1}^{\infty}\sum\limits_{a_1+a_2+a_3=i \atop a_1,a_2,a_3\geq 0}\sum\limits_{k=0}^{4a_1+2a_2}\sum\limits_{m=0}^{\infty}
\frac{(-1)^{i+a_3+k}t^{i+a_3-1}4^{a_2}}{i}\frac{i!}{a_1!a_2!a_3!}\binom{4a_1+2a_2}{k}\binom{6i+m-1}{6i-1}\\
&\qquad\qquad\qquad\qquad\qquad\qquad \times u^{a_1}v^{a_2}(u-v)^{2a_3}w^{i+a_2+2a_3+k+m}.
\end{align*}
Hence we find the coefficient of $w^n$ in $N_1^t$ is
\begin{align*}
C_1^t&=\sum\limits_{i=1}^n\sum\limits_{a_1+a_2+a_3=i \atop a_1,a_2,a_3\geq 0}\sum\limits_{k=0}^{4a_1+2a_2}
\frac{(-1)^{i+a_3+k}t^{i+a_3-1}4^{a_2}}{i}\frac{i!}{a_1!a_2!a_3!}\binom{4a_1+2a_2}{k}\\
&\qquad\qquad\qquad \times \binom{n+5i-a_2-2a_3-1-k}{6i-1}u^{a_1}v^{a_2}(u-v)^{2a_3}\\
&=\sum\limits_{i=1}^n\sum\limits_{a_1+a_2+a_3=i \atop a_1,a_2,a_3\geq 0}\sum\limits_{k=0}^{4a_1+2a_2}\sum\limits_{b=0}^{2a_3}
\frac{(-1)^{i+a_3+k+b}t^{i+a_3-1}4^{a_2}}{i}\frac{i!}{a_1!a_2!a_3!}\binom{4a_1+2a_2}{k}\\
&\qquad\qquad\qquad\times \binom{n+5i-a_2-2a_3-1-k}{6i-1}\binom{2a_3}{b}u^{a_1+2a_3-b}v^{a_2+b}\\
&=\sum\limits_{i=1}^n\sum\limits_{a_1+a_2+a_3=i \atop a_1,a_2,a_3\geq 0}\sum\limits_{b=0}^{2a_3}\frac{(-1)^{i+a_3+b}t^{i+a_3-1}4^{a_2}}{i}\frac{i!}{a_1!a_2!a_3!}\binom{n+3i-1-2a_1-a_2}{6i-1-4a_1-2a_2}\\
&\qquad\qquad\qquad\times\binom{2a_3}{b}u^{a_1+2a_3-b}v^{a_2+b}.
\end{align*}
For $N_2^t$, we have
\begin{align*}
N_2^t&=-\frac{1}{t}\sum\limits_{i=1}^{\infty}\frac{[tuw^2-t(u-v)w]^i}{i(1-w)^4}\\
&=\sum\limits_{i=1}^{\infty}\sum\limits_{m=0}^{i}\sum\limits_{a=0}^{\infty}\frac{(-1)^{i+m+1}}{i}\binom{i}{m}\binom{4i+a-1}{4i-1}t^{i-1}u^m(u-v)^{i-m}w^{i+m+a}.
\end{align*}
Hence we find the coefficient of $w^n$ in $N_2^t$ is
\begin{align*}
C_2^t&=\sum\limits_{i=1}^n\sum\limits_{m=0}^{i}\frac{(-1)^{i+m+1}}{i}\binom{i}{m}\binom{n+3i-1-m}{4i-1}t^{i-1}u^m(u-v)^{i-m}\\
&=\sum\limits_{i=1}^n\sum\limits_{m=0}^{i}\sum\limits_{b=0}^{i-m}\frac{(-1)^{i+m+b+1}}{i}\binom{i}{m}\binom{n+3i-1-m}{4i-1}\binom{i-m}{b}t^{i-1}u^{i-b}v^b.
\end{align*}
For $N_3^t$, we have
\begin{align*}
N_3^t&=-\frac{1}{t}\sum\limits_{i=1}^{\infty}\frac{[tuw^2-t(u-v)w^3]^i}{i(1-w)^4}\\
&=\sum\limits_{i=1}^{\infty}\sum\limits_{m=0}^{i}\sum\limits_{a=0}^{\infty}\frac{(-1)^{m+1}}{i}\binom{i}{m}\binom{4i+a-1}{4i-1}t^{i-1}u^{i-m}(u-v)^mw^{2i+m+a}.
\end{align*}
Hence we find the coefficient of $w^n$ in $N_3^t$ is
\begin{align*}
C_3^t&=\sum\limits_{i=1}^n\sum\limits_{m=0}^{i}\frac{(-1)^{m+1}}{i}\binom{i}{m}\binom{n+2i-1-m}{4i-1}t^{i-1}u^{i-m}(u-v)^m\\
&=\sum\limits_{i=1}^n\sum\limits_{m=0}^{i}\sum\limits_{b=0}^{m}\frac{(-1)^{m+b+1}}{i}\binom{i}{m}\binom{n+2i-1-m}{4i-1}\binom{m}{b}t^{i-1}u^{i-b}v^b.
\end{align*}
Since $H(n,3,t)=C_1^t+C_2^t+C_3^t$, we have
\begin{align*}
\sum\limits_{l,s\geq 0}W_n^3(l,s)(1-\zeta_n)^{-4s-3l}u^lv^s=-\frac{1}{n^3v}H(n,3,-1)=-\frac{1}{n^3v}(C_1^{-1}+C_2^{-1}+C_3^{-1}).
\end{align*}
We use Lemma \ref{Lemma:combinatorial identity} to show that the constant term of $H(n,3,-1)$ with respect to $v$ is zero. For a fixed non-negative integer $l$, we obtain the term $u^l$ in $H(n,3,-1)$ is
\begin{align}\label{Eq:u^l in H(n,3,-1)}
\sum\limits_{i=1}^l\frac{-1}{i}\binom{i}{l-i}\binom{n+2l-1-i}{4l-1-2i}u^l+\sum\limits_{m=0}^l\frac{(-1)^m}{l}\binom{l}{m}\binom{n+3l-1-m}{4l-1}u^l\notag\\
+\sum\limits_{m=0}^l\frac{(-1)^{m+l}}{l}\binom{l}{m}\binom{n+2l-1-m}{4l-1}u^l.
\end{align}
One can find \eqref{Eq:u^l in H(n,3,-1)} is equal to zero by using Lemma \ref{Lemma:combinatorial identity} and \eqref{Eq:zuheshu identity}.

After the concrete calculations, we find
\begin{align*}
C_1^{-1}v^{-1}&=\sum\limits_{l,s\geq 0}\sum\limits_{i=1}^n\sum\limits_{a_1+a_2=2i-s-l-1 \atop a_1,a_2\geq 0}\frac{(-1)^{s+a_2}4^{a_2}}{i}\frac{i!}{a_1!a_2!(s+l+1-i)!}\\
&\qquad\qquad \times\binom{n+2s+2l+1-i+a_2}{4s+4l+3-2i+2a_2}\binom{2s+2l+2-2i}{s+1-a_2}u^lv^s\\
&=\sum\limits_{l,s\geq 0}\sum\limits_{i=1}^{s+l+1}\frac{(-1)^s}{i}\binom{i}{s+l+1-i}\sum\limits_{a_2=0}^{2i-s-l-1}(-4)^{a_2}\binom{2i-s-l-1}{a_2}\\
&\qquad\qquad \times\binom{n+2s+2l+1-i+a_2}{4s+4l+3-2i+2a_2}\binom{2s+2l+2-2i}{s+1-a_2}u^lv^s\\
&=\sum\limits_{l,s\geq 0}\frac{(-1)^{s+1}}{s+1}\binom{s+l}{l}A(n,l,s)u^lv^s,\\
C_2^{-1}v^{-1}&=\sum\limits_{l,s\geq 0}\sum\limits_{m=0}^{l}\frac{(-1)^{s+m+1}}{s+l+1}\binom{s+l+1}{m}\binom{n+3s+3l+2-m}{4s+4l+3}\binom{s+l+1-m}{s+1}u^lv^s\\
&=\sum\limits_{l,s\geq 0}\frac{(-1)^{s+1}}{s+1}\binom{s+l}{l}\binom{n+3s+2l+2}{4s+3l+3}u^lv^s,\\
C_3^{-1}v^{-1}&=\sum\limits_{l,s\geq 0}\sum\limits_{m=0}^{l}\frac{(-1)^{s+l+m+1}}{s+l+1}\binom{s+l+1}{s+1+m}\binom{n+s+2l-m}{4s+4l+3}\binom{s+1+m}{s+1}u^lv^s\\
&=\sum\limits_{l,s\geq 0}\frac{(-1)^{s+l+1}}{s+1}\binom{s+l}{l}\binom{n+s+l}{4s+3l+3}u^lv^s.
\end{align*}
Therefore we finally get the result.
\end{proofof}

\begin{proofof}{Corollary \ref{Corollary:equivalent conj r=3}}
The proof is similar as that of Corollary \ref{Corollary:equivalent conj r=1}.
\end{proofof}


\begin{thebibliography}{99}

\bibitem{B-T-T} H. Bachmann, Y. Takeyama and K. Tasaka, Cyclotomic analogues of finite multiple zeta values, \textit{Compos. Math} \textbf{154} (2018), 2701-2721.

\bibitem{B-T-T-2} H. Bachmann, Y. Takeyama and K. Tasaka, Special values of finite multiple harmonic $q$-series at root of unity, preprint, arXiv: 1807.00411.

\bibitem{Gould} H. W. Gould, Combinatorial identities, A standardized set of tables listing 500 binomial coefficient summations, Morgantown, W. Va., 1972.  

\bibitem{Li-Pan} Z. Li and E. Pan, Sum of interpolated finite multiple harmonic $q$-series, \textit{J. Number Theory} \textbf{201} (2019), 148-175.

\bibitem{P-P-Tauraso} Kh. Pilehrood, T. Pilehrood and R. Tauraso, On $3$-$2$-$1$ values of finite multiple harmonic $q$-series at roots of unity, preprint, arXiv:2101.03576.


\end{thebibliography}
\end{document}